\documentclass[reqno]{amsart}
\newtheorem{theorem}{Theorem}[section]

\newtheorem{corollary}[theorem]{Corollary}

\theoremstyle{definition}

\newtheorem{remark}[theorem]{Remark}

\numberwithin{equation}{section}

\newcommand{\T}{\mathbb T}
\newcommand{\C}{{\mathbb C}}
\newcommand{\D}{{\mathbb D}}

\newcommand{\cS}{{\mathcal S}}

\begin{document}

\title[Boundary rigidity]{Boundary rigidity 
for some classes of meromorphic functions}

\author[V. Bolotnikov]{Vladimir Bolotnikov}
\address{Department of Mathematics \\
  College William and Mary \\
  Williamsburg, Virginia 23187-8795, U. S. A.}


\begin{abstract}
Sufficient boundary asymptotic conditions are established for a 
generalized Schur function $f$ to be identically equal to a given rational 
function $g$ unimodular on the unit circle. Similar rigidity
statements are presented for generalized Carath\'eodory
and generalized Nevanlinna functions.
\end{abstract}

\maketitle

\section{Introduction}
\setcounter{equation}{0}

In what follows, we use the following notation:
\begin{enumerate}
\item $\C$, $\D$ and $\T$ denote the complex plane, the open unit disk and
the unit circle, respectively.
\item $\cS$ --the Schur class (the closed unit ball of $H^\infty$.
\item ${\mathcal{B}_\kappa}$ --  the set of all Blaschke products of
degree $\kappa$.
\item ${\mathcal B}_p/{\mathcal B}_q$ -- the set of all coprime
quotients $g=b/\theta$ with $b\in{\mathcal B}_p$ and $\theta\in{\mathcal 
B}_q$, i.e., the set of all rational functions $g$
unimodular on $\T$ and with $p$ zeros and $q$ poles in $\D$ (counted
with multiplicities).
\item $\cS_\kappa$ -- the generalized Schur class
(introduced in \cite{kl}) consisting of all coprime quotients of the form
$f=s/b$ where $s\in\cS$ and $b\in{\mathcal B}_\kappa$.
\item $\cS_{\le \kappa}:=\bigcup_{q\le\kappa}\cS_q$ -- the set of
quotients as in (3), but not necessarily coprime.
\item $Z(f)$ -- the zero set of a function $f$.
\end{enumerate}
It is clear from definitions (4) and (7) that ${\mathcal B}_p/{\mathcal   
B}_q\subset \cS_{\le \kappa}$ whenever $q\le\kappa$.
The following rigidity result was presented in \cite{bkr} as an 
intermediate step to obtain a similar statement in the multivariable 
setting.
\begin{theorem}
Let $f\in\cS$ and let $f(z)=z+O((z-1)^4)$ as $z\to 1$. Then $f(z)\equiv z$. 
\label{T:1.1}
\end{theorem}
Generalizations and further developments can be found e.g., in 
\cite{ars, bzz, bol, chelst, elst, huang, david, tv}. Here we recall one 
from \cite{bol}.
\begin{theorem}
Let $f\in\cS$, $g\in{\mathcal B}_d$, let 
$t_1,\ldots,t_n$ be $n$ distinct points on $\T$ and let
\begin{equation}
f(z)=g(z)+o((z-t_i)^{m_i})\quad\mbox{for}\quad i=1,\ldots,n 
\label{1.2}
\end{equation}
as $z$ tends to $t_i$ nontangentially, where $m_1,\ldots,m_n$ are 
nonnegative integers. If 
\begin{equation}
\left[\frac{m_1+1}{2}\right]+\ldots +\left[\frac{m_n+1}{2}\right]>d=
\deg \, g,
\label{1.3}
\end{equation}
then $f\equiv g$. Otherwise, the uniqueness fails.
\label{T:1.2}
\end{theorem}
In \eqref{1.3}, $[x]$ denotes the largest integer that 
does not exceed a real number $x$. The last statement in Theorem 
\ref{T:1.2} means: if condition \eqref{1.3} fails for a finite  Blaschke
product $g$ and nonnegative integers $m_1,\ldots,m_n$, then for every 
choice of $n$ points $t_1,\ldots,t_n\in\T$, there are infinitely many 
functions $f\in\cS$ subject to \eqref{1.2}. Thus, conditions \eqref{1.2}
are minimal. Observe that Theorem \ref{T:1.1} follows from Theorem 
\ref{T:1.2} by letting $n=1$, $m_1=3$, $g=id$ and $t_1=1$ in the latter.

\smallskip

Since conditions \eqref{1.2} are of interpolation nature, an analog of 
Theorem \ref{T:1.2} for generalized Schur functions must exist.
It does indeed, as Theorem \ref{T:2.1} below shows. 
As we will see, this meromorphic result follows directly from its 
particular case covered by Theorem \ref{T:1.2}. Being specialized to the 
single-point case, Theorem \ref{T:2.1} gives a rigidity condition
in terms of a single asymptotic expansion (Corollary \ref{C:2.4}).
We then compare it with another single-point rigidity result recently
established in \cite{ars}. Then Caley transform will bring us to 
the analogs of Theorem \ref{T:2.1} for generalized Carath\'eodory and 
generalized Nevanlinna functions.

\section{Rigidity for generalized Schur functions}
\setcounter{equation}{0}

We start with the main result which turns out to a straightforward
consequence of Theorem \ref{T:1.2}.
\begin{theorem}
Let $\kappa,r,\ell$ be nonnegative integers and let $t_1,\ldots,t_n$ be 
$n$ distinct points on $\T$, let $g\in{\mathcal B}_\ell/{\mathcal B}_r$ 
and let us assume that a  function $f\in\cS_{\le \kappa}$ satisfies 
conditions
\begin{equation}
f(z)=g(z)+o((z-t_i)^{m_i})\quad\mbox{for}\quad i=1,\ldots,n
\label{2.4}
\end{equation}
as $z$ tends to $t_i$ nontangentially, for some 
nonnegative integers $m_1,\ldots,m_n$. If
\begin{equation}
\left[\frac{m_1+1}{2}\right]+\ldots 
+\left[\frac{m_n+1}{2}\right]>\kappa+\ell,
\label{2.5}
\end{equation}
then $f\equiv g$. 
\label{T:2.1}
\end{theorem}
{\bf Proof:} Substituting coprime quotient representations for $f$ and 
$g$
\begin{equation}  
f(z)=\frac{s_f(z)}{b_f(z)}\quad (s_f\in\cS, \; \;
b_f\in {\mathcal B}_\kappa)  
\quad\mbox{and}\quad g(z)=\frac{b(z)}{\theta(z)}\quad(b\in{\mathcal 
B}_\ell, \; \; \theta\in {\mathcal B}_r)
\label{2.6}
\end{equation}
into \eqref{2.4} and 
then multiplying both sides in \eqref{2.4} by $b_f\cdot\theta\in{\mathcal 
B}_{\kappa+r}$ we get
\begin{equation}
s_f(z)\theta(z)=b(z)b_f(z)+o((z-t_i)^{m_i})\quad\mbox{for}\quad 
i=1,\ldots,n.
\label{2.7}
\end{equation}
Since $s_f\cdot\theta\in\cS$, $b\cdot b_f\in{\mathcal B}_{\kappa+\ell}$ and 
since by \eqref{2.5}, 
${\displaystyle\sum_{i=1}^n\left[\frac{m_i+1}{2}\right]}>\kappa+\ell=\deg 
\, (b\cdot b_f)$, we conclude from \eqref{2.7} by Theorem \ref{T:1.2} that 
$s_f\cdot\theta\equiv b\cdot b_f$ 
which is equivalent, by \eqref{2.6}, to $f\equiv g$. \qed
\begin{remark} 
{\rm Observe that the membership $f\in\cS_\kappa$ means that total pole
multiplicity of $f$ does not exceed $\kappa$. Although we allow $f$
and $g$ to have different pole multiplicities, this possibility cannot 
be realized under conditions \eqref{2.5}.}   
\label{R:2.2}
\end{remark} 
Being specialized to the case $n=1$, Theorem \ref{T:2.1} gives the 
following.
\begin{corollary}
Let $\kappa, r, \ell$ be nonnegative integers, let $g\in{\mathcal 
B}_\ell/{\mathcal B}_r$ and let $f\in\cS_\kappa$ be such that 
\begin{equation}
f(z)=g(z)+o((z-t_0)^{2\kappa+2\ell+1})
\label{2.8}
\end{equation}
as $z$ tends to $t_0\in\T$ nontangentially. Then $f\equiv g$.
\label{C:2.4}
\end{corollary}
For the proof, it is enough to notice that the least integer 
$m$ satisfying inequality $\left[\frac{m+1}{2}\right]>\kappa+\ell$
is $m=2\kappa+2\ell+1$. 

\medskip

We now recall a recent result from \cite{ars} where rigidity for functions 
in $\cS_\kappa$ was established  under a 
slightly stronger condition than \eqref{2.8}.
\begin{theorem}
Let $t_0$ be a point on $\T$ and let us assume that
the numbers $\tau_0\in\T$ and $\tau_k, 
\tau_{k+1},\ldots,\tau_{2k-1}\in\C$ are such that the matrix 
${\mathbb P}=\overline{\tau}_0 TB$ is Hermitian, where 
$T$ is the lower triangular Toeplitz matrix with the bottom row
equal $[\tau_{2k-1} \; \tau_{2k-2} \; \ldots \; \tau_{k+1}, \; \tau_k]$
and $B=[b_{ij}]_{i,j=1}^k$ is the $k\times k$ right lower triangular 
matrix with the entries 
$$
b_{ij}=\left\{\begin{array}{lll}0, &\mbox{if}& 2\le i+j\le k,\\
(-1)^{j-1}{\scriptsize\left(\begin{array}{c} j-1 \\ j+i-k-1
\end{array}\right)}t_0^{j+k-1},&\mbox{if}& k+1\le i+j\le 2k.  
\end{array}\right.
$$
Let $g(z)$ be the function defined by 
\begin{equation}   
g(z)=\frac{a(z)x+b(z)}{c(z)x+d(z)}
\label{2.9}   
\end{equation}
where $x\in\T\backslash \{\tau_0\}$,
${\displaystyle
\left[\begin{array}{cc}a(z) & b(z) \\ c(z) & d(z)\end{array}\right]=
I_2-\frac{(1-z\overline{z}_0)p(z)}{(1-z\overline{t}_0)^k}
\left[\begin{array}{cc}1 & -\tau_0 \\ \overline{\tau}_0 & 
-1\end{array}\right]}$,
where $z_0\neq t_0$ is an arbitrary point on $\T$ and $p(z)$ is the 
polynomial (note that the matrix ${\mathbb P}$ is invertible by 
construction) given by $p(z)=(1-z\overline{t}_0)^k)R(z){\mathbb 
P}^{-1}R(z_0)^*$
where $R(z)=\left[\begin{array}{cccc}\frac{1}{1-z\overline{t}_0} & 
\frac{z}{(1-z\overline{t}_0)^2}& \ldots & 
\frac{z^{k-1}}{(1-z\overline{t}_0)^k}\end{array}\right]$. Then
\begin{enumerate}
\item The function $g$ is the quotient of two finite Blaschke product
with $r$ poles in $\D$ (where $r$ is
the number of negative eigenvalues of the matrix ${\mathbb P}$) and with 
the following Taylor expansion at $t_0$:
\begin{equation}   
g(z)=\tau_0+\sum_{i=k}^{2k-1}\tau_i(z-t_0)^i+O((z-t_0)^{2k}).
\label{2.10}  
\end{equation}
\item If $f\in\cS_r$ is such that 
\begin{equation}
f(z)=g(z)+O((z-t_0)^{2k+2}),  
\label{2.11}
\end{equation}
then $f\equiv g$.
\end{enumerate}
\label{T:2.5}
\end{theorem}
To embed Theorem \ref{T:2.5} into our framework we first recall 
that for every quotient of two finite Blaschke products with the 
Taylor expansion \eqref{2.10}, the matrix ${\mathbb P}$ constructed in 
the theorem is necessarily Hermitian and $\tau_0=g(t_0)$ is unimodular
(see \cite[Section 2]{bknach}. On the other hand, it follows from general 
results from \cite[Section 21]{BGR} that formula \eqref{2.9} parametrizes 
{\em all} unimodular functions $g\in{\mathcal B}_{k-r}/{\mathcal B}_r$.
Therefore, the rigidity part in Theorem \ref{T:2.5} can be 
reformulated equivalently in the following more compact form.
\begin{theorem}
Let $g\in{\mathcal B}_{k-r}/{\mathcal B}_r$ admit the Taylor expansion 
\eqref{2.10} at $t_0\in\T$. If $f\in\cS_r$
satisfies the nontangential asymptotic 
condition \eqref{2.11}, then  $f\equiv g$.
\label{T:2.6}
\end{theorem}
The main limitation in Theorem \ref{T:2.6} is that
$g$ has quite special Taylor coefficients at $t_0$ 
($\tau_1=\tau_2=\ldots=\tau_{k-1}=0$) (observe that the original 
Burns-Krantz theorem is of a different type, since there we have  
$\tau_1=1$ and $\tau_2=\tau_3=0$; however it was shown in \cite[Section 
4]{ars} that Theorem \ref{T:1.1} can be deduced from Theorem 
\ref{T:2.5}).
Corollary \ref{C:2.4} shows that rigidity holds for {\em any} quotient
of finite Blaschke products. Besides, Corollary \ref{C:2.4} shows that 
the term $O((z-t_0)^{2k+2})$ in \eqref{2.11} can be relaxed to
$o((z-t_0)^{2k+1})$, that the order of approximation can be of any parity
(not necessarily even) and that rigidity may hold also in case where only a 
bound for the pole multiplicity of $f$ is known.

\section{Rigidity for generalized Carath\'eodory and generalized 
Nevanlinna functions}
\setcounter{equation}{0}

The generalized Schur class $\cS_\kappa$ 
can be alternatively characterized as the class
of all functions $f$ meromorphic on $\D$ and such that the kernel 
$S_f(z,\zeta)=\frac{1-f(z)\overline{f(\zeta)}}
{1-z\overline{\zeta}}$ has $\kappa$ negative squares on $\D\cap{\rm 
Dom}(f)$. A related to $\cS\kappa$ is the class ${\mathcal C}_\kappa$
of generalized  Carath\'eodory functions $h$ which by definition, are 
meromorphic on $\D$ and such that the associated kernel 
$C_h(z,\zeta)=\frac{h(z)+\overline{h(\zeta)}}
{1-z\overline{\zeta}}$ has $\kappa$ negative squares on $\D\cap{\rm 
Dom}(h)$. It is convenient to include the function $h\equiv \infty$
into ${\mathcal C}_0$. Then the 
Caley transform 
\begin{equation}
f\mapsto h=\frac{1+f}{1-f}
\label{5.1}
\end{equation}
establishes a one-to-one correspondence between $\cS_\kappa$ and 
${\mathcal C}_\kappa$ and therefore, between $\cS_{\le \kappa}$ and 
${\mathcal C}_{\le \kappa}:={\displaystyle\bigcup_{r\le\kappa}{\mathcal 
C}_r}$. The representation $f=s/b$ for an $f\in\cS_\kappa$ 
combined with \eqref{5.1} implies that $h$ belongs to ${\mathcal C}_{\le 
\kappa}$ if and only if it is of the form 
\begin{equation}
h=\frac{b+s}{b-s}\quad\mbox{where}\quad b\in{\mathcal B}_\kappa, \; 
s\in\cS\quad\mbox{and}\quad Z(s)\cap Z(b)=\emptyset
\label{5.2}
\end{equation}
Theorem \ref{T:2.1} in the present setting looks as follows.
\begin{theorem}
Let $\kappa, r, \ell$ be nonnegative integers and let $g$ be of the 
form
\begin{equation}
g=\frac{b_2+b_1}{b_2-b_1}\quad\mbox{where}\quad 
b_1\in{\mathcal B}_\ell,\quad b_2\in{\mathcal B}_r
\quad\mbox{and}\quad Z(b_1)\cap Z(b_2)=\emptyset.
\label{5.3}
\end{equation}
Let us assume that a   
function $h\in {\mathcal C}_{\le\kappa}$ satisfies asymptotic equations
\begin{equation}
h(z)=g(z)+o((z-t_i)^{m_i})\quad\mbox{for}\quad i=1,\ldots,n
\label{5.4}
\end{equation}
at some points $t_1,\ldots,t_n\in\T$ and 
some nonnegative integers $m_1,\ldots,m_n$ which in turn,
are subject to \eqref{2.5}. Then $h\equiv g$.
\label{T:4.1}
\end{theorem}
{\bf Proof:} Substituting \eqref{5.2} and \eqref{5.3} into \eqref{5.4}
and then multiplying both sides in \eqref{5.4} by $(b_2-b_1)(b-s)$
we eventually get
\begin{equation}
s(z)b_2(z)=b(z)b_1(z)+o((z-t_i)^{m_i})\quad\mbox{for}\quad
i=1,\ldots,n.
\label{5.5} 
\end{equation}
Since $s\cdot b_2\in\cS$ and $b\cdot b_1\in{\mathcal B}_{\kappa+\ell}$,
we invoke Theorem \ref{T:1.2} (as in the proof of Theorem \ref{T:2.1})
to conclude from \eqref{5.5} that
$s\cdot b_2\equiv b\cdot b_1$
which implies that $h\equiv g$, thanks to \eqref{5.2} and \eqref{5.3}.\qed

\medskip

Another class related to $\cS_\kappa$ is the class ${\mathcal 
N}_\kappa$ of generalized Nevanlinna functions, that is, the functions 
$h$ meromorphic on the open upper half-plane $\C^+$ and such that the 
associated kernel 
$N_h(z,\zeta)=\frac{h(z)-\overline{h(\zeta)}}
{z-\overline{\zeta}}$ has $\kappa$ negative squares on $\C^+\cap{\rm
Dom}(h)$. The function $h\equiv\infty$ is assumed to be in ${\mathcal 
N}_0$. The classes ${\mathcal N}_\kappa$ and $\cS_\kappa$ are related 
by
\begin{equation}
h(\zeta)=i\cdot\frac{1+f(\gamma(\zeta))}{1-f(\gamma(\zeta))},\quad 
\gamma(\zeta)=\frac{\zeta-i}{\zeta+i}
\label{5.6}   
\end{equation}
which allows us to characterize ${\mathcal N}_\kappa$-functions by the 
fractional representation 
\begin{equation}
h=i\cdot \frac{b+s}{b-s}
\label{5.7}
\end{equation}  
where $s$ (analytic and bounded by one in modulus in $\C^+$) and 
$b\in{\mathcal B}_\kappa$ do not have common zeroes. For the rest of the 
paper we denote by ${\mathcal B}_k(\C^+)$ the set of finite Blaschke
products of the form 
$$
b(\zeta)=\prod_{i=1}^k\frac{\zeta-a_i}{\zeta-\bar{a}_i}\qquad 
(\zeta, \, a_i\in\C^+).
$$
Here is Theorem \ref{T:2.1} for generalized Nevanlinna functions.
\begin{theorem}
Let $\kappa, r, \ell$ be two nonnegative integers, let $g$ be of the form
\begin{equation}
g=i\cdot \frac{b_2+b_1}{b_2-b_1}\quad\mbox{where}\quad
b_1\in{\mathcal B}_\ell(\C^+),\; \;  b_2\in{\mathcal B}_r(\C^+),
\; \;  Z(b_1)\cap Z(b_2)=\emptyset.
\label{5.8}
\end{equation}
Let $\lambda_1,\ldots,\lambda_n$ be real points, let 
$m_1,\ldots,m_n$ be nonnegative integers and  let
us assume that a function $h\in {\mathcal N}_{\le\kappa}$ satisfies 
the asymptotic equations 
\begin{equation}
h(\zeta)=g(\zeta)+o((\zeta-\lambda_i)^{m_i})\quad\mbox{for}\quad 
i=2,\ldots,n
\label{5.9}
\end{equation}
as $\zeta\in\C^+$ tends to $\lambda_i$ nontangentially and
the asymptotic equation
\begin{equation}
h(\zeta)=g(\zeta)+o(|\zeta|^{-m_1}) 
\label{5.10}
\end{equation}
as $z$ tend to infinity staying inside the angle $\{z: \; \epsilon<{\rm
arg}z<\pi-\epsilon\}$. If the numbers $m_1,\ldots,m_n$ are subject to 
\eqref{2.5}, then $h\equiv g$.
\label{T:4.2}   
\end{theorem}
{\bf Proof:} Let $z:=\gamma(\zeta)$ where $\gamma$ is given in 
\eqref{5.6}. Then $t_1:=\gamma(\infty)=1\in\T$ and since 
$\lambda_i\in\T$, we have $t_i:=\gamma(\lambda_i)\in\T$ for 
$i=2,\ldots,n$. Observe that
$$
|z-t_j|=|\gamma(\zeta)-\gamma(\lambda_j)|=\frac{2|\zeta-\lambda_j|}
{|(\zeta+i)(\lambda_j+i)|}=O(|\zeta-\lambda_j|)
$$
for $j=2,\ldots,n$ and $|z-t_1|=|z-1|=|\gamma(\zeta)-1|=
\frac{2}{|\zeta+i|}=O(|\zeta|^{-1})$. Therefore, and since 
$\gamma$ maps $\C^+$ onto $\D$ conformally, we can write \eqref{5.9} and 
\eqref{5.10} as
\begin{equation}
h(\gamma^{-1}(z))=g(\gamma^{-1}(z))+o((z-t_i)^{m_i})
\quad\mbox{for}\quad i=1,\ldots,n
\label{5.11}
\end{equation}
It remains to note that the functions $-ih\circ\gamma^{-1}$ and 
$-ih\circ\gamma^{-1}$ are generalized Carathe\'eodory functions 
satisfying the assumptions of Theorem \ref{T:4.1}. Therefore, they are 
equal identically and thus, $h\equiv g$.\qed

\bibliographystyle{amsplain}

\end{document}